\documentclass{amsart}[10pt]
\usepackage{enumerate, amsfonts,epsfig,amssymb, amsmath, amsthm}

\newlength{\tabwidth}
\newlength{\tabheight}
\newlength{\tabrule}
\newlength{\tabwidthx}
\newlength{\tabheightx}

\def\gentabbox#1#2#3#4{\vbox to \tabheight{\setlength{\tabrule}{#3}%
  \setlength{\tabwidthx}{#1\tabwidth}\addtolength{\tabwidthx}{\tabrule}%

\setlength{\tabheightx}{#2\tabheight}\addtolength{\tabheightx}{-\tabheight}%
  \hbox to #1\tabwidth{%
 \hspace{-0.5\tabrule}\rule{\tabrule}{#2\tabheight}\hspace{-\tabrule}%
    \vbox to #2\tabheight{\hsize=\tabwidthx%
      \vspace{-0.5\tabrule}\hrule width\tabwidthx height\tabrule%
      \vspace{-0.5\tabrule}\vfil%
      \hbox to \tabwidthx{\hss#4\hss}%
        \vfil\vspace{-0.5\tabrule}%
      \hrule width\tabwidthx height\tabrule\vspace{-0.5\tabrule}}%
 \hspace{-\tabrule}\rule{\tabrule}{#2\tabheight}\hspace{-0.5\tabrule}}%
  \vspace{-\tabheightx}}}
\def\genblankbox#1#2{\vbox to \tabheight{\vfil\hbox to
#1\tabwidth{\hfil}}}
\def\tabbox#1#2#3{\gentabbox{#1}{#2}{0.4pt}{\strut #3}}

\catcode`\:=13 \catcode`\.=13 \catcode`\;=13
\catcode`\>=13 \catcode`\^=13
\def:#1\\{\hbox{$#1$}}
\def.#1{\tabbox{1}{1}{$#1$}}
\def>#1{\tabbox{2}{1}{$#1$}}
\def^#1{\tabbox{1}{2}{$#1$}}
\def;{\genblankbox{1}{1}\relax}
\catcode`\:=12 \catcode`\.=12 \catcode`\;=12
\catcode`\>=12 \catcode`\^=7

\newenvironment{tableau}{\bgroup\catcode`\:=13 \catcode`\.=13
  \catcode`\;=13 \catcode`\>=13 \catcode`\^=13
  \setlength{\tabheight}{3ex}\setlength{\tabwidth}{3ex}%
  \def\b##1##2##3{\gentabbox{##1}{##2}{1.2pt}{\vbox{##3}}}%
  \def\n##1##2##3{\gentabbox{##1}{##2}{0.4pt}{\vbox{##3}}}%
  \vbox\bgroup\offinterlineskip}{\egroup\egroup}

\newtheorem{theorem}{Theorem}[section]
\newtheorem{corollary}[theorem]{Corollary}
\newtheorem{lemma}[theorem]{Lemma}
\newtheorem{proposition}[theorem]{Proposition}
\newtheorem{conjecture}[theorem]{Conjecture}
\newtheorem*{theorem*}{}

\theoremstyle{definition}
\newtheorem{definition}[theorem]{Definition}

\theoremstyle{remark}

\newtheorem{fact}{Fact}[section]
\newtheorem{example}[theorem]{Example}
\newtheorem{remark}[theorem]{Remark}

\begin{document}
\title[Cells and Constructible Representations]{Cells and Constructible Representations\\ in type $B$}
\author{Thomas Pietraho}
\email{tpietrah@bowdoin.edu} 
\address{Department of Mathematics\\Bowdoin College\\Brunswick,
Maine 04011}
 \maketitle

\begin{abstract} We examine the partition of a finite Coxeter group of
type $B$ into cells determined by a weight function $L$. The main
objective of these notes is to reconcile Lusztig's description of
constructible representations in this setting with conjectured
combinatorial descriptions of cells.
\end{abstract}
\section{Introduction}
Consider a finite Coxeter group $W$ together with a weight function
$L:W\rightarrow \mathbb{Z}$, as in \cite{lusztig:unequal}. Every
weight function is specified by its values on the simple reflections
in $W$ and defines an  Iwahori-Hecke algebra $\mathcal{H}$ by
explicit generators and relations.  Furthermore, following Lusztig,
a weight function determines a partition $W$ into left, right, and two-sided cells,
each one of which carries representations of $\mathcal{H}$ and $W$
\cite{lusztig:unequal}.  Their role in
the representation theory of reductive algebraic groups over finite
or $p$-adic fields is described in Chapter 0 of
\cite{lusztig:unequal}.  Cells also arise in the study of
rational Cherednik algebras and the Calogero-Moser space, see
\cite{gordon:calogero} and \cite{gordon:quiver}.

Left cell representations of $W$ are intimately related to its
constructible representations; that is, the minimal class of
representations of $W$ which contains the trivial representation and
is closed under truncated induction and tensoring with sign.  In
fact, left-cell and constructible representations coincide when $L$
is the length function on $W$, see \cite{lusztig:cells}. With
the additional stipulation that the conjectures (P1)-(P15) of
\cite{lusztig:unequal} hold, M.~Geck has shown this to be true for
general weight functions as well \cite{geck:constructible}.

Left cells are well understood for dihedral groups and Coxeter
groups of type $F_4$. We focus our attention on the remaining case of  Coxeter groups of type $B_n$. The weight function is then specified by two integer parameters $a$ and $b$:

\begin{center}
\begin{picture}(300,30)
\put( 50 ,10){\circle*{5}} \put( 50, 8){\line(1,0){40}} \put( 50,
12){\line(1,0){40}} \put( 48, 20){$b$} \put( 90 ,10){\circle*{5}}
\put(130 ,10){\circle*{5}} \put(230 ,10){\circle*{5}} \put( 90,
10){\line(1,0){40}} \put(130, 10){\line(1,0){25}} \put(170,
10){\circle*{2}} \put(180, 10){\circle*{2}} \put(190,
10){\circle*{2}} \put(205, 10){\line(1,0){25}} \put( 88 ,20){$a$}
\put(128, 20){$a$} \put(222, 20){$a$}
\end{picture}
\end{center}
Given $a,b \neq 0$, we may assume both are positive by
\cite{ginzburg:cherednik}(5.4.1), and write $s=\tfrac{b}{a}$ for
their quotient. Parameterizations of the left, right, and two-sided
cells of $W$ have been obtained by Garfinkle  \cite{garfinkle3} in
the equal parameter case $s=1$, by Lusztig \cite{lusztig:leftcells}
and Bonnaf\'e, Geck, Iancu, and Lam \cite{bgil} for $s=\tfrac{1}{2}$
and $s=\frac{3}{2}$, and Bonnaf\'e--Iancu \cite{bonnafe:iancu} and
Bonnaf\'e \cite{bonnafe:two-sided} in the asymptotic case $s>n-1$.
Furthermore, a description for the remaining values of $s$ has been
conjectured by Bonnaf\'e, Geck, Iancu, and Lam in \cite{bgil}. On
the other hand,  constructible representations of $W$ were already described
by Lusztig for all values of $s$ by relying on conjectures
(P1)-(P15) of \cite{lusztig:unequal}.

The above parametrizations of cells in type $B_n$ can be stated in terms of families of standard domino tableaux of arbitrary rank. Reconciling
this description of cells with Lusztig's parametrization of
constructible representations is therefore a natural question and is
the main purpose of this paper.  We focus our attention on the case
$s \in \mathbb{N}$, excluding the cases for which the left cell and
constructible representations are conjectured to be irreducible.
In this setting, our main result shows the consistency of Lusztig's conjectures (P1)-(P15) with the conjectural descriptions of cells.  As corollaries to this work, we amend the original conjectural description of two-sided Kazhdan-Lusztig cells in \cite{bgil}, and examine under what circumstances Lusztig's notion of special representation can exist in the unequal parameter case.

The paper is organized as follows.  Section \ref{section:definitions} defines cells in unequal parameter Hecke algebras and summarizes the requisite combinatorics.  In Section \ref{section:cells}, we examine the conjectural combinatorial description of cells in Weyl groups of type $B_n$ and its consequences.  Finally, Section \ref{section:constructible} connects this work with constructible representations.

\section{Definitions and Preliminaries}
\label{section:definitions}
We begin by defining Kazhdan-Lusztig cells for unequal parameter Hecke algebras.  The rest of the section is devoted to the combinatorics pertinent to their conjectured parametrization in type $B_n$.  Our main goal is to describe certain properties of cycles in a domino tableau, first on the level of partitions, and then on the level of symbols.

\subsection{Kazhdan-Lusztig Cells} Let $(W,S)$ be a Coxeter system with a
weight function $L:W \rightarrow \mathbb{Z}$ which takes positive
values on all $s \in S$. Define $\mathcal{H}$ to be the generic
Iwahori-Hecke algebra over $\mathcal{A}= \mathbb{Z}[v, v^{-1}]$ with
parameters $\{v_s \, | \, s \in S\}$, where  $v_w = v^{L(w)}$ for
all $w \in W$. The algebra $\mathcal{H}$ is free over $\mathcal{A}$
and has a basis $\{T_w \, | \, w \in W\}$ in terms of which
multiplication takes the form
$$T_s T_w = \left\{
        \begin{array}{ll}
            T_{sw} & \text{if $\ell(sw) > \ell(w)$, and}\\
            T_{sw}+(v_s-v_s^{-1}) T_w & \text{if $\ell(sw) < \ell(w)$.}
        \end{array}
        \right.
        $$
for $s \in S$ and $w \in W$. As in \cite{lusztig:unequal}(5.2), it
is possible to construct a Kazhdan-Lusztig basis of $\mathcal{H}$
which we denote by $\{C_w \; | \; w \in W\}$.  For $x,y \in W$ and
some $h_{xyz} \in \mathcal{A}$, multiplication in $\mathcal{H}$
takes the form
$$ C_x C_y = \sum_{z \in W} h_{xyz} C_z.$$

\begin{definition}(\cite{lusztig:unequal}(8.1)) Fix $(W,S)$  a Coxeter system with a weight
function $L$.

\begin{enumerate}

    \item We will say $w' \leq_{\mathcal{L}} w$ if there exists $s \in S$ for
which $C_{w'}$ appears with a non-zero coefficient in $C_s C_w$ and
reuse the same notation $\leq_{\mathcal{L}}$ for the transitive
closure of this binary relation. The equivalence relation associated with
the preorder $\leq_{\mathcal{L}}$ will be denoted by $\sim_\mathcal{L}$ and its
equivalence classes will be called  {\it Kazhdan-Lusztig left cells} of $W$.

    \item We will say $w' \leq_{\mathcal{R}} w$ iff $w'^{-1}
\leq_{\mathcal{L}} w^{-1}$, write $\sim_\mathcal{R}$ for the
corresponding equivalence relation and call its
equivalence classes the {\it Kazhdan-Lusztig right cells} of $W$.

    \item Finally, we define $\leq_{\mathcal{LR}}$ as the pre-order generated
by $\leq_{\mathcal{L}}$ and $\leq_{\mathcal{R}}$, write
$\sim_\mathcal{LR}$ for the corresponding equivalence relation and
call its equivalence classes the {\it Kazhdan-Lusztig two-sided cells}
of $W$.
\end{enumerate}
\end{definition}

Each Kazhdan-Lusztig cell carries a representation of the
Iwahori-Hecke algebra $\mathcal{H}$. We reconstruct the definition
of \cite{lusztig:unequal}(8.3). If $\mathfrak{C}$ is a
Kazhdan-Lusztig left cell and $w \in \mathfrak{C}$, then
$$[\mathfrak{C}]_\mathcal{A} = \bigoplus_{w' \leq_\mathcal{L} w} \mathcal{A} C_{w'} \Big/ \bigoplus_{w' \leq_\mathcal{L} w, w' \notin \mathfrak{C} } \mathcal{A} C_{w'},$$
is a quotient of two left ideals in $\mathcal{H}$ and therefore is a
left $\mathcal{H}$-module.  The set $[\mathfrak{C}]_\mathcal{A}$
does not depend on the specific choice of $w \in \mathfrak{C}$, is
free over $\mathcal{A}$, and has a basis $\{e_w \; | \; w \in
\mathfrak{C}\}$ indexed by elements of $\mathfrak{C}$ where $e_w$ is
the image of $C_w$ in the above quotient.  Elements of $\mathcal{H}$
act on $[\mathfrak{C}]_\mathcal{A}$ via $C_x e_y = \sum_{z \in
\mathfrak{C}} h_{xyz} e_z$ for $x \in W$ and $y \in \mathfrak{C}$.
Finally, by  restricting to scalars, $[\mathfrak{C}]_\mathcal{A}$
gives rise to a $W$-module which we denote by $[\mathfrak{C}]$.  The
situation is similar for right cells.  If $\mathfrak{D}$ is a
Kazhdan-Lusztig two-sided cell and $w \in \mathfrak{D}$, then
$$[\mathfrak{D}]_\mathcal{A} = \bigoplus_{w' \leq_\mathcal{LR} w} \mathcal{A} C_{w'} \Big/ \bigoplus_{w' \leq_\mathcal{LR} w, w' \notin \mathfrak{D} } \mathcal{A} C_{w'},$$
is a quotient of two two-sided ideals of $\mathcal{H}$ and therefore
is a  $\mathcal{H}$-bimodule.  The set $[\mathfrak{D}]_\mathcal{A}$
does not depend on the specific choice of $w \in \mathfrak{D}$, is
free over $\mathcal{A}$, and has a basis $\{e_w \; | \; w \in
\mathfrak{D}\}$ indexed by elements of $\mathfrak{D}$ where $e_w$ is
the image of $C_w$ in the above quotient. Again by  restricting to
scalars, $[\mathfrak{D}]_\mathcal{A}$ gives rise to a $W$-module
which we denote by $[\mathfrak{D}]$.

\subsection{Partitions}

Let $p = (p_1, p_2, \ldots, p_k)$ be a partition of $m$ with the
convention  $p_1 \geq p_2 \geq \ldots \geq p_k >0$.   We will
routinely identify a partition $p$ with its Young diagram $Y_p$, or
a left-justified array of boxes whose lengths decrease from top to
bottom.  Thus the partition $(4,3,3,1) = (4,3^2,1)$ will correspond
to the Young diagram

$$
\begin{tiny}
\begin{tableau}
:.{}.{}.{}.{}\\
:.{}.{}.{}\\
:.{}.{}.{}\\
:.{}\\
\end{tableau}
\end{tiny}
$$

If the Young diagram of a partition can be tiled by dominos, we will
say that the underlying partition is of {\it rank zero}.  In
general, suppose that we can remove a domino from a Young diagram in
such a way that what remains is another Young diagram justified at
the same row and column.  Repeating this process starting with a
partition $p$ will eventually terminate, and the reminder will be a
Young diagram of a partition $(r, r-1, r-2 , \ldots, 1)$ for some $r
\geq 0$.  We will write $p \in \mathcal{P}_r$ and say that $p$ is of
rank $r$.  The rank of a partition is unique; each partition of $m$
belongs to $\mathcal{P}_r$ for exactly one value of $r$.  The {\it
core} of $p$ is the triangular partition $(r, r-1, r-2 , \ldots,
1).$

Let $s_{ij}$ denote the square in row $i$ and column $j$ of a Young
diagram and extend this notion somewhat by letting $i$ and $j$ take
values that may not describe squares in the Young diagram. We define
two sets of squares related to the Young diagram of a partition $p$.

\begin{definition}  For a partition $p \in \mathcal{P}_r$ and its corresponding
Young diagram $Y_p$, consider squares $s_{ij}$ such that $i+j \equiv
r \; (\textup{mod } 2)$, $i+j>r+1$, and either the addition of $s_{ij}$ to $Y_p$
or the removal of $s_{ij}$ from $Y_p$ yields another Young diagram.
Among these, we will say
    \begin{itemize}
        \item  $s_{ij} \in \mathcal{C}(p)$ iff $i$ is odd, and
        \item $s_{ij} \in \mathcal{H}(p)$ iff $i$ is even.
    \end{itemize}
We will write $\mathcal{HC}(p)$ for the union $\mathcal{C}(p) \cup
\mathcal{H}(p)$.  Furthermore, we define  sets $\mathcal{HC}^*(p),$
$\mathcal{C}^*(p),$ and $\mathcal{H}^*(p)$ exactly as above, but
requiring $i+j>r+2$ instead.  We will say $s_{ij} \in
\mathcal{HC}(p)$ is {\it filled} if it lies in $Y_p$ itself, otherwise, we
will say it is {\it empty}.  Finally, let $\gamma_p = | \{s_{ij} \in Y_p \, | \, i+j=r+2\}|$ and $\kappa_p$ be the number of filled squares in $\mathcal{HC}(p).$
\end{definition}

\begin{fact} \label{fact:corners} We have the following easy consequences:
    \begin{enumerate}
        \item The element of $\mathcal{HC}(p)$ with the smallest row number lies in $\mathcal{C}(p).$
        \item Elements of $\mathcal{C}(p)$ and $\mathcal{H}(p)$ alternate with increasing row number.
        \item $|\mathcal{H}(p)| \leq |\mathcal{C}(p)|.$
        \item Both $\gamma_p = r+1 $ as well as $\kappa_p \neq 0$ imply that $\mathcal{HC}(p)=\mathcal{HC}^*(p).$
    \end{enumerate}
\end{fact}

\begin{example}  In the Young diagram of the rank 2 partition $p=(4,3^2,1)$,
these sets are $\mathcal{C}(p) =\mathcal{C}^*(p)=
\{s_{15},s_{33},s_{51}\}$ and $\mathcal{H}(p)
=\mathcal{H}^*(p)=\{s_{24},s_{42}\}$.
\end{example}

\begin{definition}  Let $s_{ij}, s_{kl}, \text{ and $s_{mn}$}  \in \mathcal{HC}(p)$.  We will say that $s_{mn}$ lies between $s_{ij}$ and $s_{kl}$ iff
$m$ is between $i$ and $k$ and $n$ is between $j$ and $l$ (where $m$
is between $i$ and $k$ iff $i \leq m \leq k$ or $i \geq m \geq k$).  We will say that two squares of $\mathcal{HC}(p)$ are {\it adjacent  in } $\mathcal{HC}(p)$ if no other square of $\mathcal{HC}(p)$ lies between them.
\end{definition}

\begin{remark}\label{remark:corners}
The  set $\mathcal{HC}(p)$ coincides with the union of the sets of corners and holes of $p$ as defined for $r=0$ and $1$ in \cite{garfinkle1}.  When $r=0$, $\mathcal{C}(p)$  is the set of corners and $\mathcal{H}(p)$ is the set of holes; however, our definitions diverge from Garfinkle's for $r=1.$  This is not unexpected: when $r=1$, the left cells identified in Conjecture \ref{conjecture:bgil} are not the ones studied by Garfinkle.  Their parametrization depends on a different choice of variable squares,  see Remark \ref{remark:cycles}.
\end{remark}

\begin{remark}
The set $\mathcal{HC}^*(p)$ is
precisely the union of {\it addable} and {\it removable} squares as
defined in \cite{gordon:calogero}. Following \cite{gordon:calogero},
the {\it heart} of $p$ will be the partition obtained from $p$ by
removing all filled squares of $\mathcal{HC}^*(p)$.
\end{remark}

\subsection{Tableaux}
Consider a partition $p$.  A {\it domino tableau of shape $p$} is a Young diagram of $p$ whose
squares are labeled by a set $M$ of non-negative integers in such a way that every positive integer labels exactly two adjacent squares and all labels increase weakly along both rows and columns.  A domino tableau is {\it standard} if $M \cap \mathbb{N} = \{1, \ldots, n\}$ for some $n$,  and of {\it rank $r$} if $p \in \mathcal{P}_r$ and $0$ labels the square $s_{ij}$ iff $i+j < r+2.$  We will write $SDT_r(p)$ for the family of all standard domino
tableaux of rank $r$ and shape $p$  and $SDT_r(n)$ for the family of
all standard domino tableaux of rank $r$ which contain exactly $n$ dominos.

The moving-though operation on a domino tableau defines
another domino tableau whose labels agree on a certain
subset of its squares.  In a domino tableau of rank $r$, we will say
that the square $s_{ij}$ in row $i$ and column $j$ is {\it variable}
iff $i+j \equiv r \; (\textup{mod }2)$; otherwise, we will say it is {\it
fixed}.

\begin{remark}
\label{remark:cycles} Our choice of fixed and variable squares coincides with that of \cite{garfinkle1} when $r=0$ but not when $r=1$.  As mentioned in Remark \ref{remark:corners}, in this latter case the left cells identified in Conjecture \ref{conjecture:bgil} are not the ones studied by Garfinkle and their description relies on the above assignment of fixed and variable squares.
\end{remark}

 Consider a domino $D=D(k,T)$ with label $k \in \mathbb{N}$ in a
domino tableau $T$ of rank $r$ and write $supp \, D(k,T)$ for the set of its
underlying squares. Then $supp \, D(k,T)$ contains a fixed and a variable
square. Suppose that we wanted to create another domino
tableau by changing the label of the variable square of $D$ in a way
that preserved the labels of all fixed squares of $T$ while
perturbing $T$ minimally. As in \cite{garfinkle1}, this leads to the notion of a cycle in a
domino tableau; we define it presently.

\begin{definition}
Suppose that  $supp \, D(k,T)= \{s_{ij},s_{i+1,j}\}$ or
$\{s_{i,j-1},s_{ij}\}$ and the square  $s_{ij}$ is fixed. Define
$D'(k)$ to be a domino labeled by the integer $k$ with $supp \,
D'(k,T)$ equal to
\begin{align*}
        \{s_{ij}, s_{i-1,j}\}     & \text{ if $k< label \, s_{i-1,j+1}$,
        and} \\
        \{s_{ij}, s_{i,j+1}\}    & \text{ if $k> label \,
        s_{i-1,j+1}$.}
\end{align*}
Alternately, suppose that  $supp \, D(k,T)= \{s_{ij},s_{i-1,j}\}$ or
$\{s_{i,j+1},s_{ij}\}$ and the square  $s_{ij}$ is fixed. Define
$supp \, D'(k,T)$ to be
\begin{align*}
      \{s_{ij},s_{i,j-1}\}    & \text{ if $k< label \, s_{i+1,j-1}$,
      and} \\
         \{s_{ij},s_{i+1,j}\}        & \text{ if $k> label \,
         s_{i+1,j-1}.$}
\end{align*}
\end{definition}

\begin{definition}
The cycle $c=c(k,T)$ through $k$ in a domino tableau $T$ of rank $r$ is
a union of labels of dominos in $T$  defined by the condition that
$l \in c$ if either $l=k$, or either $supp \, D(l,T) \cap supp \,
D'(m,T) \neq
        \emptyset$ or $supp \, D'(l,T) \cap supp \,
D(m,T) \neq
        \emptyset$ for some $D(m,T) \in c$.
\end{definition}

We will refer to the set of dominos with labels in a cycle $c$ as
the cycle $c$ itself.   For a domino tableau $T$ of rank
$r$ and a cycle $c$ in $T$, define $MT(T,c)$ by
replacing every domino $D(l,T) \in c$ by the corresponding domino
$D'(l,T)$.  It follows that $MT(T,c)$ is domino tableau, and in general, the shape of
$MT(T,c)$ will either equal the shape of $T$, or one square will be
removed (or added to the core) and one will be added
\cite{garfinkle1}(1.5.27). A cycle $c$ is called {\it closed} in the
former case and {\it open} in the latter.   We will write $OC(T)$ for the
set of open cycles in $T$.  For $c\in OC(T)$, we will write $S_b(c)$
for the square that is either removed from the shape of $T$ or added
to the core of $T$ by moving through $c$. Similarly, we will write
$S_f(c)$ for the square that is added to the shape of $T$. Note that
$S_b(c)$ and $S_f(c)$ are always variable squares.  Consistent with
Garfinkle's notation in \cite{garfinkle3}, we will write $OC^*(T)$
for the set of non-core open cycles in $T$, that is, cycles for
which both $S_b(c)$ and $S_f(c)$ lie in $\mathcal{HC}^*(p)$ with
$p=shape \, T$.  For a cycle in $OC^*(T)$, $S_b(c) \in
\mathcal{C}^*(p)$ and $S_f(c) \in \mathcal{H}^*(p)$, or $S_b(c) \in
\mathcal{H}^*(p)$ and $S_f(c) \in \mathcal{C}^*(p)$.

Let $\mathcal{C} $ be a set of cycles in a domino tableau $T$ of rank $r$. According
to \cite{garfinkle1}(1.5.29), moving through disjoint
 cycles in a domino tableau are independent operations, allowing us to unambiguously write
$MT(T,\mathcal{C})$ for the domino tableau obtained by simultaneously moving-through all
of the cycles in the set $\mathcal{C}$.  If $\mathcal{C} \subset OC^*(T)$, then $MT(T,\mathcal{C})$ is another domino tableau of rank $r$, $\mathcal{C} \subset OC^*(MT(T,\mathcal{C}))$, and  $MT(MT(T,\mathcal{C}),\mathcal{C})=T.$ If $\mathcal{C} = OC(T) \setminus OC^*(T),$ then $MT(T,\mathcal{C})$ can be interpreted as a domino tableau of rank $r+1$ and comes endowed with new sets of fixed and variable squares and consequently, cycles.

\begin{definition} For a standard domino tableau of rank $r$, we define the cycle structure set of $T$
 as the set of ordered pairs $cs(T)$
consisting of the beginning and final squares of every cycle in $T$
and $cs^*(T)$ as the restriction of this set to non-core open
cycles. That is:
\begin{align*}
cs(T) & = \{ (S_b(c), S_f(c)) \; | \; c \in OC(T) \} \text{, and } \\
cs^*(T)&  = \{ (S_b(c), S_f(c)) \; | \; c \in OC^*(T) \}.
\end{align*}
Finally, write $\widetilde{cs}(T)$ and $\widetilde{cs}^*(T)$ for the sets obtained from the above by changing their underlying ordered pairs into unordered pairs.
\end{definition}

We would like a similar notion for partitions that does not directly rely on an underlying tableau.  First note that if $p = shape \, T$, then $\mathcal{HC}(p)$ consists exactly of the $\kappa_p$ beginning and $\kappa_p$ final squares of non-core open cycles of $T$, the $\gamma_p$ final squares of core open cycles of $T$, and the $r+1-\gamma_p$ empty squares adjacent to the core of $Y_p$; consequently, we have $|\mathcal{HC}(p)|=2\kappa_p+r+1$.

\begin{definition} Consider $p \in \mathcal{P}_r$.
A {\it cycle structure set} $\sigma$ for $p$, or alternately, for $\mathcal{HC}(p)$, is a pairing of squares in $\mathcal{H}^*(p)$ with
squares in $\mathcal{C}^*(p)$ for which \label{definition:csset}
    \begin{enumerate}
        \item exactly $\gamma_p$ squares remain unpaired, and \label{definition:csset:one}
        \item every square $c \in \mathcal{HC}^*(p)$ which lies between
        $a$ and $b$ for a pair $\{a,b\} \in \sigma$ must be paired with
        another square which lies between $a$ and $b$.
        \label{definition:csset2}
    \end{enumerate}
\end{definition}

\begin{example} The partition $(4,3^2,1)$ of rank $r=2$ admits exactly four cycle structure
 sets:  $\{\{s_{15},s_{24}\}\},$ $\{\{s_{24},s_{33}\}\},$
 $\{\{s_{33},s_{42}\}\},$ and $\{\{s_{42},s_{51}\}\}.$
\end{example}

Note that a  cycle structure set for $p$ contains exactly $\kappa_p$ pairs. Cycle structure sets for tableaux and partitions are closely related. Given a standard domino tableau $T$ of rank $r$, the set $\sigma=\widetilde{cs}^*(T)$ is a cycle structure set for the partition $p=shape(T) \in \mathcal{P}_r$ by elementary properties of open cycles.  Conversely, as detailed in the following proposition, a cycle structure set for an arbitrary partition $p$ always arises as a cycle structure set
for some domino tableau.

\begin{proposition}  If $p \in \mathcal{P}_r$ and $\sigma$
is a cycle structure set for $p$, then there exists a standard domino tableau $T$ of rank $r$
with  $\widetilde{cs}^*(T)=\sigma$. \label{proposition:exists}
\end{proposition}
\begin{proof}

We proceed by induction on the number of pairs in $\sigma$.  If
$\sigma$ is empty, then any $T \in SDT_r(p)$ suffices.  Otherwise, consider a pair
$\{s,s'\} \in \sigma$. By Definition
\ref{definition:csset}(\ref{definition:csset2}), $s$ and $s'$ can be
chosen in such a way that they are adjacent in $\mathcal{HC}(p)$.
Lemma 3.4 of \cite{gordon:calogero} implies that there exists a
standard domino tableau $T'$ of rank $r$ and a non-core open cycle
$c$ in $T'$ with $S_b(c)=s$ and $S_f(c)=s'$ whose dominos form a rim
ribbon $R$ of $T'$.  If $R$ contains $t$ dominos, then by the proof
of \cite{gordon:calogero}(3.4), $T'$ can be chosen in such a way that the
dominos in $T'\setminus R$ are labeled by elements of the set $\{1,
2, \ldots, n-t\}.$  Now by Definition
\ref{definition:csset}(\ref{definition:csset2}), $\sigma \setminus
\{s, s'\}$ is either empty or contains a pair $\{s'',s'''\}$ adjacent in $\mathcal{HC} (shape \, (T'\setminus R))$ and the proposition  follows by induction.
\end{proof}

\begin{proposition}  If $p \in \mathcal{P}_r$ and $\mathcal{S}$ is a subset of
$\mathcal{HC}^*(p)$ consisting of $\kappa_p$ elements, then there exists a cycle structure set $\sigma$ for $p$ where each pair in $\sigma$ contains exactly one element of
$\mathcal{S}.$ \label{proposition:cs2}
\end{proposition}

\begin{proof}
When $|\mathcal{S}|=\kappa_p=0$, the only cycle structure set for $p$ is empty.  Assuming that $|\mathcal{S}|=\kappa_p \geq 1$, Fact \ref{fact:corners} implies $\mathcal{HC}(p)=\mathcal{HC}^*(p)$.  Therefore, we can find a pair $\{s,t\}$ of squares adjacent in $\mathcal{HC}(p)$ with $s\in \mathcal{S}$ and $t\notin \mathcal{S}$.  Again by Fact \ref{fact:corners}, $\mathcal{C}(p)$ and $\mathcal{H}(p)$ alternate with increasing row number and consequently one of $s$ and $t$ must lie in $\mathcal{C}^*(p)$ and the other in $\mathcal{H}^*(p).$  Working recursively, this pair can be extended to a pairing of elements of $\mathcal{H}^*(p)$ with elements of $\mathcal{C}^*(p)$ which satisfies the properties of a cycle structure set.
\end{proof}

\subsection{Symbols}\label{section:symbols}

A symbol of defect $s$ is an array of numbers of the form
$$
\Lambda= \left(
\begin{array}{cccccc}
  \lambda_1 & \lambda_2  &  &\ldots  & & \lambda_{N+s} \\
  {} & \mu_1 & \mu_2 & \ldots  & \mu_N & {}
\end{array}
\right)
$$
where $\{\lambda_i\}$ and $\{\mu_i\}$ are, perhaps empty, strictly increasing
sequences of non-negative integers.  Define an equivalence relation
on the set of symbols of defect $s$ by letting $\Lambda$ be
equivalent to the symbol
$$
\Lambda'= \left(
\begin{array}{ccccccc}
  0 & \lambda_1+1  & \lambda_2+2 &\ldots  & & \lambda_{N+s} +N+s\\
   {} & 0 & \mu_1+1 & \ldots  & \mu_N+N & {}
\end{array}
\right)
$$
We will write $Sym_s$ for the set of equivalence classes of symbols
of defect $s$. It is possible to define a map from partitions to
symbols via the following procedure.  Given a partition $p = (p_1,
p_2, \ldots, p_k)$, form an extended partition  $p^\sharp=(p_1, p_2,
\ldots, p_{k'})$ by adding an additional zero term to $p$ if the
rank of $p$ has the same parity as $k$.  The set
$\{p_i+k'-i\}_{i=1}^{k'}$ can be divided into odd and even parts
$\{2\mu_i+1\}_{i=1}^N$ and $\{2\lambda_i\}_{i=1}^{N+s}$ from which
the symbol $\Lambda_p$ corresponding to $p$ can be constructed by
arranging the $\lambda_i$ and $\mu_i$ as above.  We will write
$\widetilde{p}_i$ for the entry of $\Lambda_p$ determined from the
part $p_i$ of $p^\sharp$.

Let $\mathcal{P}^2$ be the set of ordered pairs of partitions, and
write $\mathcal{P}^2(n)$ for the subset of $\mathcal{P}^2$ where the
sum of parts of both partitions sum to $n$. Given a symbol of defect
$s$, it is also possible to construct an ordered pair of partitions.
With $\Lambda$ as above, let $d_\Lambda = \{\lambda_i - i
+1\}_{i=1}^{N+s}$ and $f_\Lambda =\{\mu_i - i +1\}_{i=1}^N$. The
following is an immediate consequence of \cite{james:kerber}(2.7).

\begin{theorem}
The maps $p \mapsto \Lambda_p$ and  $\Lambda \mapsto (d_{\Lambda},
f_{\Lambda})$ define  bijections
$$\mathcal{P}_r \rightarrow Sym_{r+1} \rightarrow
\mathcal{P}^2$$ for all values of $r$.  Furthermore, the composition
of these two maps yields a bijection between $\mathcal{P}_r(n)$ and
$\mathcal{P}^2(n)$. \label{theorem:bijections}
\end{theorem}

Since the set $\mathcal{P}^2(n)$ parameterizes the irreducible
representations of the Weyl group $W_n$ of type $B_n$, the above can
be used to identify irreducible $W_n$-modules with symbols of fixed
defect, as in \cite{lusztig:unequal}, or with partitions of fixed
rank. It is the latter interpretation that we employ in the
following sections. We will write $[p]$ and $[\Lambda]$ for the
irreducible $W_n$-modules associated to the partition $p$ and the symbol
$\Lambda$ in this manner.

We would like to understand the correspondence of Theorem
\ref{theorem:bijections} in slightly greater detail. For a symbol
$\Lambda$, we write $Z_1(\Lambda)$ for the set of entries that
appear once among its rows, and $Z_2(\Lambda)$ for the set of
entries that appear twice.

\begin{lemma}
The set $Z_1(\Lambda_p)$ of single entries of the symbol $\Lambda_p$
can be identified with the parts of $p^\sharp$ whose rows in the
Young diagram $Y_p$ end in a, perhaps empty, square of
$\mathcal{HC}(p)$. This establishes a bijective map $Z_1(\Lambda_p)
\leftrightarrow \mathcal{HC}(p).$ \label{lemma:z1chp}
\end{lemma}
\begin{proof}
We first show that elements of $Z_2(\Lambda_p)$ arise  from the rows of
$Y_p$ which do not terminate in a square of $\mathcal{HC}(p)$. Note
that $z \in Z_2(\Lambda_p)$ implies $z=\widetilde{p}_i=\widetilde{p}_{i+1}$ for some $i$, and furthermore $p_i$ must equal $p_{i+1}$. A parity argument shows that if the row of $p_i$ in $Y_p$ ends in a fixed square, $\widetilde{p}_i$ will differ from $\widetilde{p}_{i+1}$. Hence every $z\in Z_2(\Lambda_p)$ must correspond to a pair of consecutive equal parts $p_i,p_{i+1}$ of $p^\sharp$ where the row of $p_i$ in $Y_p$ ends in a variable square. It is easy to check that $\widetilde{p}_i=\widetilde{p}_{i+1}$ for such a pair.  Hence
elements of $Z_2(\Lambda_p)$ correspond to pairs of consecutive equal rows of
$Y_p$, the first of which ends in a variable square. Such pairs yield
precisely the rows of $Y_p$ which do not terminate in a square of
$\mathcal{HC}(p)$ and the lemma follows.
\end{proof}

\section{Combinatorial Cells}
\label{section:cells}

This section examines equivalence relations on the Weyl group of
type $B_n$ defined via a Robinson-Schensted algorithm and standard
domino tableaux.  As stated more explicitly in the next section, the equivalence classes which they define are expected to coincide with unequal parameter Kazhdan-Lusztig cells in type $B_n$.

\subsection{Robinson-Schensted Algorithms}

The Weyl group $W_n$ of type $B_n$ is the group of permutations of
the set $\{\pm 1, \pm 2, \ldots , \pm n \}$ which commute with the
involution $i \mapsto -i$.   Generalized Robinson-Schensted maps
$G_r : W_n \rightarrow SDT_r(n) \times SDT_r(n)$ defined in
\cite{garfinkle1} and \cite{vanleeuwen:rank} construct bijections
between elements of the Weyl group of type $B_n$ and same-shape
pairs of standard domino tableaux of rank $r$ for each non-negative
integer $r$.  We will write $G_r(w) = (S_r(w), T_r(w))$ for the
image of an element $w$ and refer to the components of the ordered
pair as the rank $r$ left and right tableaux of $w$.

There is a natural description of the relationship between the
bijections $G_r$ for differing $r$ described in terms of the moving
though map for open cycles, see \cite{pietraho:rscore}. We also
point out that for $r$ sufficiently large, $G_r$ recovers another
generalization of the Robinson-Schensted algorithm for
hyperoctahedral groups defined in \cite{stanton:white} and
\cite{okada}.  See \cite{pietraho:equivalence} for a more detailed
description.

\subsection{Combinatorial Left, Right, and Two-Sided Cells}

\begin{definition} Consider $x,y \in W_n$ of type $B_n$ and fix a non-negative integer
$r$.  We will say
\begin{enumerate}
    \item $x \approx_{\mathcal{L}} y$ if their right tableaux of rank $r$ are
related by moving through some set of non-core open cycles, that is,
iff $T_r(x) = MT(T_r(y), \mathcal{C})$ for some $\mathcal{C} \subset
OC^*(T_r(y))$,

    \item  $x \approx_{\mathcal{R}} y$ iff $x^{-1}
\approx_{\mathcal{L}} y^{-1}$, and

    \item write $x \approx_{\mathcal{LR}} y$ for the relation generated by
$\approx_\mathcal{L}$ and $\approx_\mathcal{R}$.

\end{enumerate}
\end{definition}

\begin{center}
\begin{figure}
\epsfig{file=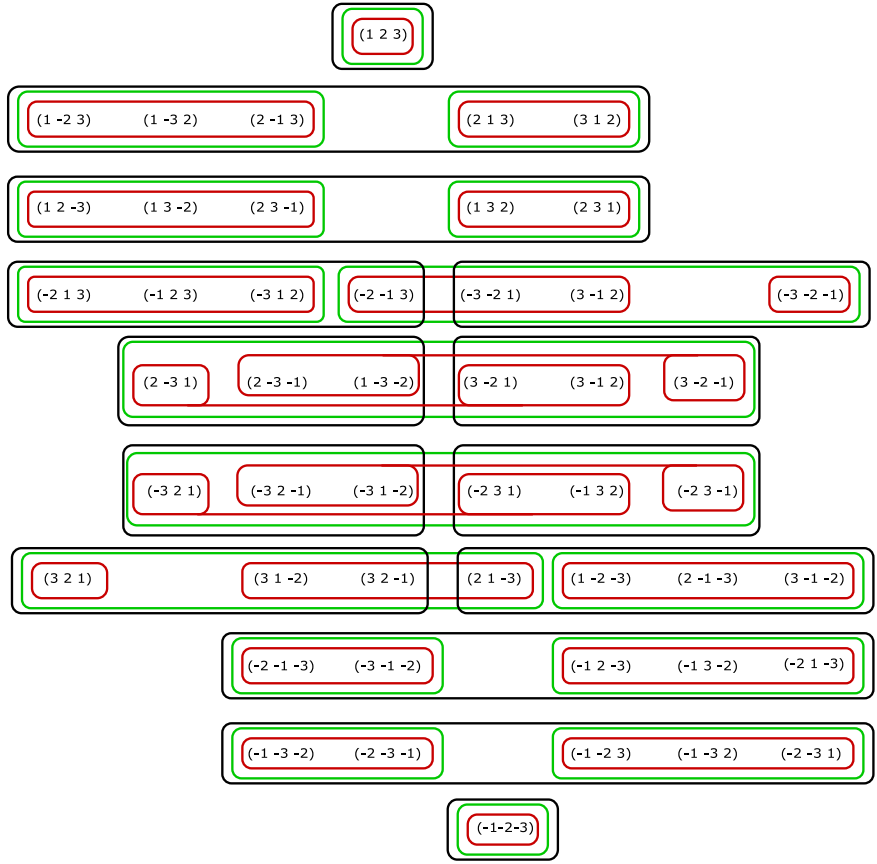, height=4.5in} \caption{Combinatorial
left cells in $W_3$.  Black represents cells for rank $r=0$, green
represents $r=1$, and red represents $r\geq 2$. The cells are not
successive refinements for increasing values of the partition rank
parameter $r$.}\label{figure:leftcells}
\end{figure}
\end{center}

\begin{center}
\begin{figure}
\epsfig{file=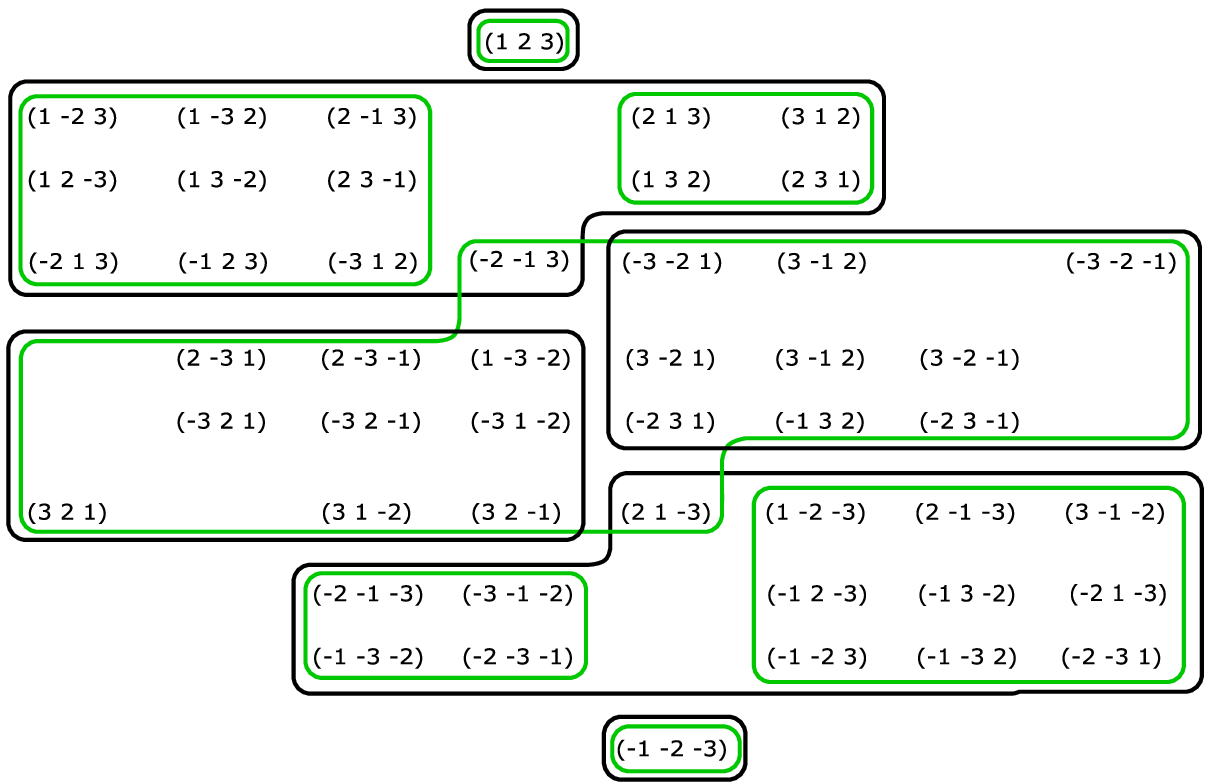,height=3in}
\epsfig{file=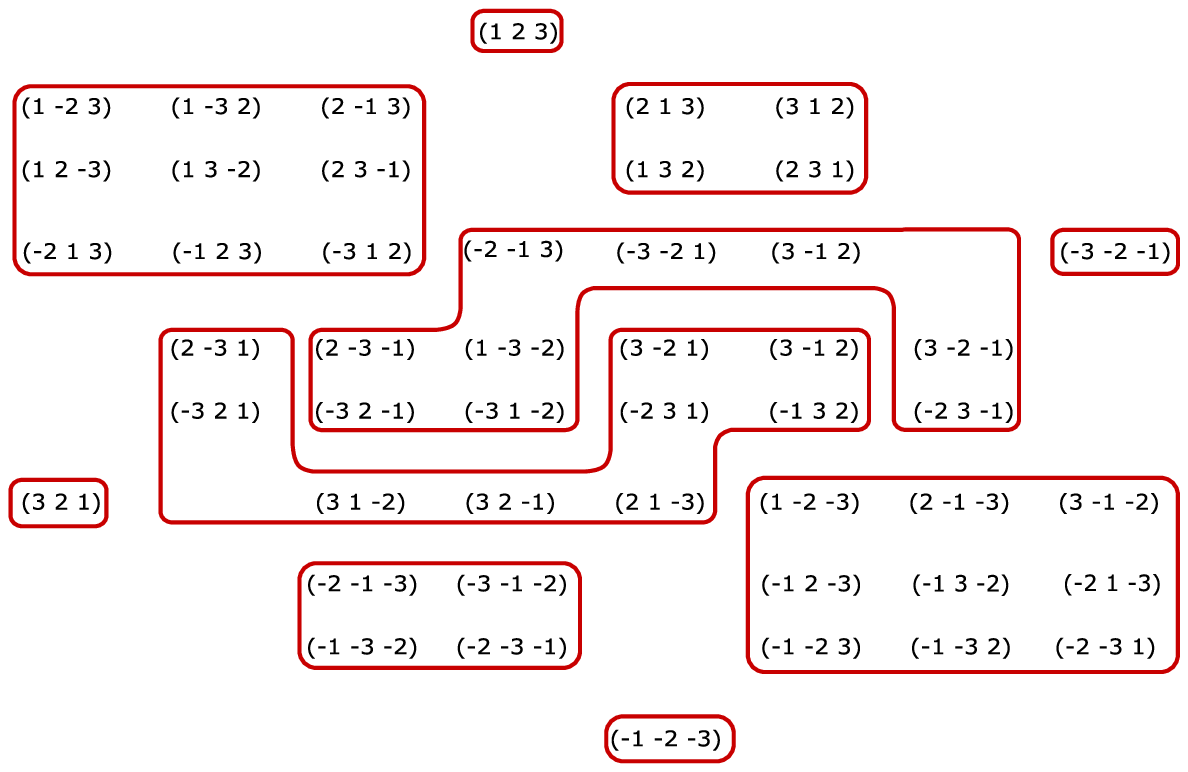,height=3in} \caption{Combinatorial
double cells in $W_3$.  Black represents cells for rank $r=0$, green
represents $r=1$, and red represents $r\geq 2$. The cells are not
successive refinements for increasing values of the partition rank
parameter $r$.} \label{figure:twosidedcells}
\end{figure}
\end{center}

Defined in this way, we will call the equivalence classes of $\approx_\mathcal{L}$,
$\approx_\mathcal{R}$, and $\approx_\mathcal{LR}$ in $W_n$ {\it reducible
combinatorial left, right, and two-sided cells of rank $r$}.  Within this paper, we will generally omit the adjective ``reducible" of this definition.  Although we
suppress it in the notation, the cells depend on the choice of
parameter $r$.   By \cite{vanleeuwen:rank}(4.2), the map $w \mapsto
w^{-1}$ on $W_n$ carries combinatorial left cells to combinatorial
right cells and preserves combinatorial two-sided cells. As seen in
Figures \ref{figure:leftcells} and \ref{figure:twosidedcells},
combinatorial cells do not behave simply with respect to a change in
$r$, although it is possible to describe a precise relationship
\cite{pietraho:rscore}. When $r>n-2$, the situation is somewhat
simpler.  There are no non-core open cycles, implying both, that
combinatorial left cells are determined simply by right tableaux,
and by the main result of \cite{pietraho:rscore}, that for these values of $r$, all
combinatorial cells are actually independent of $r$.

\vspace{-.45in}


The main result of \cite{pietraho:equivalence} shows that combinatorial left cells admit the
following alternate description.  A similar characterization holds
for combinatorial right cells.
\begin{theorem}[\cite{pietraho:equivalence}]  Combinatorial left
cells in the Weyl group of type $B_n$  are generated by the
equivalence relations of having the same right tableau in either
rank $r$ or rank $r+1$.
\end{theorem}

\begin{proposition}
Consider $x,y \in W_n$, fix a non-negative integer $r$, and let $p$
and $p'$ be the  shapes of $T_r(x)$ and $T_r(y)$ respectively. Then
$x \approx_\mathcal{LR} y$ iff $\mathcal{HC}(p)=\mathcal{HC}(p')$
and $p$ and $p'$ differ only in the choice of filled squares in
$\mathcal{HC}(p)$. \label{proposition:twosided}
\end{proposition}
\begin{proof}
Suppose $x \approx_\mathcal{LR} y$.  Note that if either $x'
\approx_\mathcal{L} x''$ or $x' \approx_\mathcal{R} x''$, then
$T_r(x')$ and $T_r(x'')$ or $S_r(x')$ and $S_r(x'')$ differ by
moving through a, perhaps empty, set of non-core open cycles.  Since
moving through non-core open cycles acts on the level of partitions
by only changing which squares are filled in $\mathcal{HC}(p)$, the
forward direction of the above follows.

For the other direction, first note that two elements whose tableaux
are of the same shape are necessarily in the same combinatorial
two-sided cell:  if $T$ and $S$ are two standard dominos of the same
shape, then
$$G_r^{-1}(X,T) \approx_{\mathcal{L}} G_r^{-1}(T,T)
\approx_{\mathcal{R}} G_r^{-1}(S,T) \approx_{\mathcal{L}}
G_r^{-1}(S,S) \approx_{\mathcal{R}} G_r^{-1}(Y,S)$$ for all $X$ and
$Y$ of the same shape.  The rest of the proof follows as in
\cite{gordon:calogero}(3.5).
\end{proof}

\subsection{Tableau shapes of elements within combinatorial cells}

We examine more closely the sets of partitions that appear among
shapes of tableaux of elements in combinatorial cells. Fix a
combinatorial left cell $\mathfrak{C}$ and a combinatorial two-sided
cell $\mathfrak{D}$ of rank $r$.  Write $\pi(\mathfrak{C})$ and
$\pi(\mathfrak{D})$ for the sets of partitions that appear among
tableaux shapes of their elements.

It is clear that the rank $r$ right tableaux of the elements of $\mathfrak{C}$ share a common cycle structure set and $\pi(\mathfrak{C})$ consists exactly of those partitions derived from a choice of a filled square in each of its constituent pairs.  If $k_\mathfrak{C}$ is the number of non-core open cycles in the tableaux of the elements of $\mathfrak{C}$, then $|\pi(\mathfrak{C})| = 2^{k_\mathfrak{C}}$. The partitions in $\pi(\mathfrak{D})$ can be determined via the following observations.
\begin{enumerate}
        \item According to
            Proposition \ref{proposition:twosided}, the sets
            $\mathcal{HC}(p)$, $\mathcal{H}(p)$, and $\mathcal{C}(p)$ are
            constant among
            $p \in \pi(\mathfrak{D})$. We will emphasize this by writing
            $\mathcal{HC}_\mathfrak{D}$, $\mathcal{H}_\mathfrak{D}$, and
            $\mathcal{C}_\mathfrak{D}$ for these sets.
        \item The number of filled squares in $\mathcal{HC}_\mathfrak{D}$ is constant on $\pi(\mathfrak{D})$.  We will denote this number by $\kappa_\mathfrak{D}$.  It is equal to $\kappa_p$ for any $p\in \pi(\mathfrak{D})$.
    \end{enumerate}
From Proposition \ref{proposition:twosided}, partitions in $\pi(\mathfrak{D})$ are determined by choices of $\kappa_\mathfrak{D}$ filled squares among $\mathcal{HC}_\mathfrak{D}$. Consequently,  $\pi(\mathfrak{D})$ contains exactly ${|\mathcal{HC}_\mathfrak{D}| \choose k_\mathfrak{D}}$ partitions.
The next proposition points out a fundamental difference in the relationship between the
sets $\pi(\mathfrak{C})$ and $\pi(\mathfrak{D})$ in the two cases
when $r=0$ or $r>n-2$, and when $0<r\leq n-2$.

\begin{proposition}
Consider a combinatorial two-sided cell $\mathfrak{D}$. The
intersection $$I_\mathfrak{D}=\bigcap_{\mathfrak{C} \subset
\mathfrak{D}} \pi(\mathfrak{C})$$ is non-empty iff $\kappa_\mathfrak{D} =
0 $ or $\kappa_\mathfrak{D}= |\mathcal{H}_\mathfrak{D}|$, in which case
it contains a unique partition. In particular, this occurs for all
combinatorial two-sided cells of rank $r=0$ and $r>n-2$.
\label{proposition:distinguished}
\end{proposition}
\begin{proof}

We first claim that if $p \in I_\mathfrak{D}$, then all filled squares of $\mathcal{HC}(p)$ must lie in $\mathcal{H}^*(p)$.  This is trivially true if $\kappa_p=0$, so assume otherwise and note that by Fact \ref{fact:corners},  $\mathcal{HC}(p)=\mathcal{HC}^*(p)$. Suppose $c \in \mathcal{C}(p)=\mathcal{C}^*(p)$. Since $|\mathcal{HC}(p)|= 2 \kappa_p + r+1$ and $|\mathcal{H}(p)|\leq |\mathcal{C}(p)|,$ we can choose a set of $\kappa_p$ elements of $\mathcal{C}^*(p)$ which excludes $c$.  By Proposition \ref{proposition:cs2}, there is a cycle structure set for $p$ which leaves $c$ unpaired and by Proposition \ref{proposition:exists}, there is a standard domino tableau $T$ with $\widetilde{cs}^*(T) = \sigma$.  If we let $\mathfrak{C}$ be the combinatorial left cell associated to $T$, then no partition in $\pi(\mathfrak{C})$ has the square $c$ filled, and consequently all filled squares of $\mathcal{HC}(p)$ must lie in $\mathcal{H}^*(p)$.

Armed with this observation, suppose first that  $\kappa_\mathfrak{D} = |\mathcal{H}_\mathfrak{D}| \neq 0$ so that $\mathcal{HC}_\mathfrak{D}= \mathcal{HC}^*_\mathfrak{D}$.  Every cycle
structure set on $\mathcal{HC}_\mathfrak{D}$ must pair all of the $\kappa_\mathfrak{D}$ squares within $\mathcal{H}_\mathfrak{D}$, implying that the partition with precisely all squares of $\mathcal{H}_\mathfrak{D}$ filled lies in $\pi(\mathfrak{C})$ for all cells $\mathfrak{C} \subset \mathfrak{D}$. Furthermore, this partition is unique in $I_\mathfrak{D}$:  any partition in $\pi(\mathfrak{D})$ must have exactly $\kappa_\mathfrak{D}$ filled squares in $\mathcal{HC}_\mathfrak{D}$ and, as observed above, for partitions in $I_\mathfrak{D}$ these must lie in $\mathcal{H}_\mathfrak{D}$.

When $k_\mathfrak{D} <|\mathcal{H}_\mathfrak{D}|$, Proposition
\ref{proposition:cs2} can be used to construct a cycle structure set
$\sigma$ on $\mathcal{HC}_\mathfrak{D}$ which leaves an arbitrary
$h \in \mathcal{H}_\mathfrak{D}$ unpaired.  We can associate a
combinatorial left cell $\mathfrak{C}$ to $\sigma$ as above and note
that $h$ is empty in every partition of $\pi(\mathfrak{C})$.  Since
$h$ was arbitrary, any partition appearing in $I_\mathfrak{D}$ must
have all $h \in \mathcal{H}_\mathfrak{D}$ empty.  Hence unless $\kappa_\mathfrak{D}=0$, $I_\mathfrak{D}$ is empty. When $k_\mathfrak{D}=0$,  $\pi(\mathfrak{D})$ consists of a unique partition and $|I_\mathfrak{D}|=1$. Finally, we note that if $r=0$, $k_\mathfrak{D} = |\mathcal{H}_\mathfrak{D}|$ and if $r>n-2$, $k_\mathfrak{D}=0$.
\end{proof}

\begin{remark}\label{remark:special}
In the case $r=0$, the unique partition in $I_\mathfrak{D}$ is
called special and corresponds to Lusztig's notion of special
representation of $W_n$ under the map defined by
Theorem \ref{theorem:bijections}, see \cite{lusztig:leftcells}.  A
consequence of the above proposition is that similarly distinguished
partitions do not exist for a range of values of $r$.  When
interpreted in terms of the conjectures describing the
Kazhdan-Lusztig cells in type $B_n$ stated in the next section, and combined with the results of \cite{pietraho:module}, this precludes the existence of distinguished representations of $W_n$ in the general unequal parameter case.
\end{remark}

\begin{example}
Consider the partition $p=(4,3^2,1)$ of rank $2$.  Elements of $W_4$
whose tableaux have shape $p$ lie in a combinatorial two-sided cell
$\mathfrak{D}$. The set
$\mathcal{HC}_\mathfrak{D}=\mathcal{HC}^*_\mathfrak{D}$ equals
$\{s_{15},s_{24},s_{33}, s_{42},s_{51}\}$, with only the square
$s_{33}$  filled, hence $\kappa_\mathfrak{D}=1$. Consequently, listing
the partitions of $\pi(\mathfrak{D})$ entails deciding which square
of $\mathcal{HC}_\mathfrak{D}$ is filled.  The possible partitions
are $(5,3,2,1)$, $(4^2,2,1)$, $(4,3^2,1)$, $(4,3,2^2)$, and
$(4,3,2,1^2)$.

The shapes of elements in combinatorial left cells contained in this
combinatorial two-sided cell fall into the following four
categories: $\{(5,3,2,1),(4^2,2,1)\},$ $\{(4^2,2,1), (4,3^2,1)\}, $
$\{(4,3^2,1), (4,3,2^2)\},$ and $\{(4,3,2^2)\}, (4,3,2,1^2)\},$ each
corresponding to a choice of a cycle structure set on
$\mathcal{HC}_\mathfrak{D}$. In particular, it is clear that no
partition is common to all of these sets.
\end{example}

\begin{remark}For every combinatorial left cell $\mathfrak{C}$,  $\pi(\mathfrak{C})$ admits a natural structure of an
elementary abelian $2$-group. Since the right tableau of any element
in $\mathfrak{C}\subset\mathfrak{D}$ is of the form
$MT(T,\mathcal{C})$ for some $C \subset OC^*(T)$,
$\pi(\mathfrak{C})$ is determined entirely by the positions of each
of the $\kappa_\mathfrak{D}$ non-core open cycles of $T$, and corresponds
to the choices of a filled square within each pair of
$\widetilde{cs}^*(T)$. Because the moving-through operations on
cycles in $T$ are independent,  a choice of a distinguished
partition in $\pi(\mathfrak{C})$ defines a natural structure of an
elementary abelian 2-group of order $2^{\kappa_\mathfrak{D}}$.  For
$r=0$, this is described in \cite{mcgovern:leftcells}.
\end{remark}

\section{Kazhdan-Lusztig cells and constructible representations}
\label{section:constructible}

We examine the relationship of Kazhdan-Lusztig cells and
combinatorial cells in type $B_n$ and reconcile
Lusztig's description of constructible representations with
combinatorial cells.  We restrict our attention to the case where
the defining parameter $s$ is an integer, focusing on the case when the
conjectured cells and constructible representations are not
irreducible.  The key will be the results of Section \ref{section:definitions} relating partitions and symbols.

\subsection{Cells in type $B_n$} \label{subsection:cellsbn}

We restrict the setting to the Weyl group of type $B_n$ with
generators as in the following diagram:
\begin{center}
\begin{picture}(300,30)
\put( 50 ,10){\circle*{5}} \put( 50, 8){\line(1,0){40}} \put( 50,
12){\line(1,0){40}} \put( 48, 20){$t$} \put( 90 ,10){\circle*{5}}
\put(130 ,10){\circle*{5}} \put(230 ,10){\circle*{5}} \put( 90,
10){\line(1,0){40}} \put(130, 10){\line(1,0){25}} \put(170,
10){\circle*{2}} \put(180, 10){\circle*{2}} \put(190,
10){\circle*{2}} \put(205, 10){\line(1,0){25}} \put( 88 ,20){$s_1$}
\put(128, 20){$s_2$} \put(222, 20){$s_{n-1}$}
\end{picture}
\end{center}
Suppose the weight function $L$ is defined by $L(t)=b$ and
$L(s_i)=a$ for all $i$.  We will examine the case when $\frac{b}{a}
\in \mathbb{N}$, and set $s=\frac{b}{a}$.  The following is a
conjecture of Bonnaf\'e, Geck, Iancu, and Lam, and appears as
Conjecture B in \cite{bgil}:

\begin{conjecture}[\cite{bgil}]\label{conjecture:bgil} Consider a Weyl group of type $B_n$
with a weight function $L$ and parameter $s$ defined as above.
Kazhdan-Lusztig left, right, and two-sided cells for parameter $s\in
\mathbb{N}$ coincide with combinatorial left, right, and two-sided
cells of rank $s-1.$
\end{conjecture}

This conjecture is well-known to be true for $s=1$ by work of
Garfinkle \cite{garfinkle3}, and has been verified when $s>n-1$ by
Bonnaf\'e and Iancu, \cite{bonnafe:iancu} and Bonnaf\'e
\cite{bonnafe:two-sided}. It has also been shown to hold for all
values of $s$ when $n \leq 6$, see \cite{bgil}.  The above is
restated more explicitly as Conjecture D in \cite{bgil}.  It implicitly assumes the existence of a partition which is distinguished in the sense of Remark \ref{remark:special} within the partitions arising among tableaux of elements of a two-sided cell.  However,
in light of Proposition \ref{proposition:distinguished} such a partition does not exist in general and 
the characterization of Kazhdan-Lusztig two-sided cells in Conjecture D must be rephrased using the description of combinatorial two-sided cells of Proposition \ref{proposition:twosided}.

\subsection{Constructible Representations}

The set of constructible representations $Con(W)$ of a Weyl group $W$ is the
smallest class of representations which contains the trivial
representation and is closed under truncated induction and tensoring
with the sign representation, see \cite{lusztig:unequal}(22.1).
Although this is not clear from the notation, this set depends on
the weight function chosen to define $\mathcal{H}$.  For the results
of this section, we assume that Lusztig's conjectures P1-P15 of
\cite{lusztig:unequal}(14.2) are true.  Under this assumption, M.~Geck has described  the relationship between constructible representations and
Kazhdan-Lusztig left cells:

\begin{proposition}[\cite{geck:constructible}] Consider a finite Coxeter group $W$ with a weight function $L$ and let $\mathfrak{C}$ be a Kazhdan-Lusztig left cell in $W$ defined from $L$.  Then
  \begin{enumerate}
        \item $[\mathfrak{C}]$ is a constructible $W$-module, and
        \item every constructible $W$-module can be obtained in this
    way.
    \end{enumerate}
\end{proposition}

Let us again restrict the setting to the Weyl group of type $B_n$
with weight function $L$ defining a parameter $s$. We begin our
description of constructible representations by first recalling the
one of Lusztig \cite{lusztig:unequal}(22.6). Let $\Lambda$ be a
symbol of defect $s$ and let $Z_1=Z_1(\Lambda)$ and
$Z_2=Z_2(\Lambda)$.  If $Y\subset Z_1$, define a new symbol
$$\Lambda_Y = \left(
\begin{array}{c}
Z_2  \sqcup Z_1 \setminus Y \\
Z_2  \sqcup  Y
\end{array}
\right)
$$
We would like to restrict the set of subsets $Y$ for which this
construction will be carried out. An involution $\iota: Z_1
\rightarrow Z_1$ is  {\it admissible} iff
\label{definition:admissible}

    \begin{enumerate}
        \item it contains exactly $s$ fixed points,
        \item whenever  $z'\in Z_1$ lies strictly between
$z$ and $\iota(z)$ for any $z \in Z_1$, then $z'$ is not a fixed
point and $\iota(z')$ lies strictly between $z$ and
$z'$.\label{definition:admissible2}

    \end{enumerate}
Given an admissible involution $\iota$, define a set $S_\iota$
consisting of subsets of $Z_1$ by letting $Y \in S_\iota$ iff it
contains exactly one element from each orbit of $\iota$. Recalling
the parametrization of $W_n$-modules by symbols of defect $s$ from
Section \ref{section:symbols}, define a $W_n$-module by
$$c(\Lambda,\iota) = \bigoplus_{Y \in S_\iota} [\Lambda_Y]$$
The modules $c(\Lambda,\iota)$ and $c(\Lambda',\iota')$ are equal
iff $\Lambda$ and $\Lambda'$ have the set of entries and also
$\iota=\iota'$ .

\begin{proposition}[\cite{lusztig:unequal}(22.23)]
Consider a symbol $\Lambda$ together with an admissible involution $\iota$.
Then
    \begin{enumerate}
        \item $c(\Lambda,\iota)$ is a constructible $W_n$-module, and
        \item every constructible $W_n$-module can be obtained in this
    way.
    \end{enumerate}
\end{proposition}

Now consider a partition $p \in \mathcal{P}_{s-1}$.  If $Y$ is a
subset of $\mathcal{HC}(p)$, let $p_Y$ be the partition obtained
from the heart of $p$ by filling exactly the squares of
$\mathcal{HC}(p)$ which correspond to $Y$. Given a cycle structure
set $\sigma$ for $p$, define a set $\mathcal{S}_\sigma$ consisting of
subsets of $\mathcal{HC}(p)$ by letting $Y \in \mathcal{S}_\sigma$
iff $Y$ contains exactly one element from each pair in $\sigma$.
Recalling the parametrization of $W_n$-modules by partitions of rank
$s-1$ from Section \ref{section:symbols}, we define a $W_n$-module by
$$c(p,\sigma) = \bigoplus_{Y \in S_\sigma} [p_Y]$$
The modules $c(p,\sigma)$ and $c(p',\sigma')$ are equal iff $p$ and
$p'$ have the same heart and $\sigma=\sigma'$ .  The $W_n$-modules
obtained in this way are precisely the constructible ones.

\begin{theorem}Consider a partition $p \in \mathcal{P}_{s-1}$ and a cycle structure set
$\sigma$ for $p$, then
\begin{enumerate}
    \item $c(p,\sigma)$ is a constructible $W_n$-module, and
    \item every constructible $W_n$-module can be obtained in this
    way.
\end{enumerate}
\end{theorem}

\begin{proof}
Construct $\Lambda_p$, a symbol of rank $s$, as in Theorem
\ref{theorem:bijections}. We first show that
$c(p,\sigma)=c(\Lambda_p, \iota)$ for some admissible involution
$\iota$.  When $\mathcal{HC}(p)\supsetneq\mathcal{HC}^*(p)$, we have
$|\mathcal{HC}(p)|=s$ and Definition
\ref{definition:csset}(\ref{definition:csset:one}) implies that the
only cycle structure set $\sigma$ on $p$ is trivial.  Hence
$c(p,\sigma)=[p]$.  By Lemma \ref{lemma:z1chp}, the corresponding
symbol $\Lambda_p$ will have $|Z_1(\Lambda_p)| = s$, implying that
the only admissible involution $\iota$ on $Z_1$ is trivial. Hence
$c(\Lambda_p,\iota)=[\Lambda_p]=[p]=c(p,\sigma).$

Thus we assume $\mathcal{HC}(p)=\mathcal{HC}^*(p)$ and write
$Z_1=Z_1(\Lambda_p)$ and $Z_2=Z_2(\Lambda_p)$. We describe a
bijection between the cycle structure sets for $p$ and admissible
involutions $\iota:Z_1\rightarrow Z_1$.  Images of the orbits of
$\iota$ under the map of Lemma \ref{lemma:z1chp} form a pairing
$\sigma_\iota$ on the squares of $\mathcal{HC}(p)$. Noting that the
squares in $\mathcal{HC}(p)$ alternate between $\mathcal{H}(p)$ and
$\mathcal{C}(p)$ with increasing row number, Definition
\ref{definition:admissible}(\ref{definition:admissible2}) implies
that $\sigma_\iota$ is in fact a pairing between squares of
$\mathcal{H}(p)$ and $\mathcal{C}(p)$. Furthermore, it follows
directly from the definition that $\sigma_\iota$ is in fact a cycle
structure set for $p$. This process is easily reversed, establishing
the desired bijection. Write $\iota_\sigma$ for the admissible
involution associated with the cycle structure set $\sigma$.

We would like to show that $c(p,\sigma)=c(\Lambda_p, \iota_\sigma).$
Lemma \ref{lemma:z1chp} establishes a bijection between
$S_{\iota_\sigma}$ and $S_\sigma$.   If $\widetilde{Y}$ represents
the image of $Y \in S_{\iota_\sigma}$, it is sufficient to show that
the symbol  $(\Lambda_p)_Y = \Lambda_{p_{\widetilde{Y}}}$ for all $Y
\in S_{\iota_\sigma}$.  It is clear that
$Z_1((\Lambda_p)_Y)=Z_1(\Lambda_{p_{\widetilde{Y}}})$ and
$Z_2((\Lambda_p)_Y)=Z_2(\Lambda_{p_{\widetilde{Y}}}).$ Consider a
square $s_{ij} \in \widetilde{Y}$ and write $\iota_{ij}$ for the
corresponding element of $Y \subset Z_1$.  It is enough to show show
that $\iota_{ij}$ appears in the bottom row of the symbol
$\Lambda_{p_{\widetilde{Y}}}$.  With $k'$ defined as in Section
\ref{section:symbols}, note that $j+k'-i$ is odd.  By the definition
of the map $p \rightarrow \Lambda_p$ and Lemma \ref{lemma:z1chp},
$\iota_{ij}$ must equal $\tfrac{j+k'-i-1}{2}$ and hence appears in
the bottom row of $\Lambda_{p_{\widetilde{Y}}}$, as desired.

Finally, since the map of Theorem \ref{theorem:bijections} is a
bijection and we've established a bijection between cycle structure
sets and $s$-admissible involutions,  every constructible $W_n$-module
appears as $c(p,\sigma)$ for some $p$ and $\sigma$, since it appears
as $c(\Lambda,\iota)$ for some $\Lambda$ and $\iota$.

\end{proof}

The above theorem can easily be restated in terms of tableaux. To
each tableau $T \in SDT_r(n)$ we associate a $W_n$-module $[T]$ in the
following manner. For each family of open cycles $\mathcal{C}$ in
$T$ define $p_\mathcal{C}$ to be the shape of the tableau
$MT(T,\mathcal{C})$ obtained from $T$ by moving through
$\mathcal{C}$. Let
$$
[T]= \bigoplus_{\mathcal{C} \subset OC^*(T)} [p_\mathcal{C}]
$$
The partitions which can be obtained by moving through non-core open
cycles in a tableau depend only on the cycle structure of the
tableau, hence the modules $[T]$ and $[T']$ are equal iff the
underlying partitions have the same heart and $\widetilde{cs}^*(T) =
\widetilde{cs}^*(T')$.  The $W_n$-modules obtained in this way are
precisely the constructible ones.

\begin{corollary}  \label{corollary:module:tableau}Consider a standard domino tableau $T$, then
\begin{enumerate}
    \item $[T]$ is a constructible $W_n$-module, and
    \item every constructible $W_n$-module can be obtained in this
    way.
\end{enumerate}
\end{corollary}

\begin{proof}  The module $[T]$ is precisely $c(shape(T),\widetilde{cs}^*(T))$,
and hence constructible.  Conversely, a constructible module
$c(p,\sigma)$ equals $[T]$ for some tableau $T$ by Proposition
\ref{proposition:exists}.
\end{proof}

Given a Coxeter system $(W,S)$ with weight function $L$, a {\it
family} of partitions is an equivalence class defined by the
transitive closure of the relation linking $p$ and $p'$  iff $[p]$
and $[p']$  appear as simple components of some constructible
representation of $W$. The following relates families and the
partitions appearing in combinatorial two-sided cells.  It is a
version of the result of \cite{geck:constructible}(4.3).

\begin{proposition}
Consider $W$ is of type $B_n$ with a weight function $L$ and
parameter $s$ and let $\mathfrak{D}$ be a combinatorial two-sided
cell.   Then the family of $p\in \pi(\mathfrak{D})$ is precisely
$\pi(\mathfrak{D})$.
\end{proposition}

\begin{proof}
If $p,p'$ lie in the same family, then they must have the same
heart, implying $p,p' \in \pi(\mathfrak{D})$. We show the converse.
If $k_\mathfrak{D}=0$, then $|\pi(\mathfrak{D})|=1$, and we are
done. Otherwise, note that $\mathcal{HC}(p) = \mathcal{HC}^*(p)$ and we can let $p^\Uparrow$ be the partition with the same
heart as $p$ but with the top-most $k_\mathfrak{D}$ squares of
$\mathcal{HC}(p)$ filled. We will show that $p$ and
$p^\Uparrow={p'}^\Uparrow$ lie in the same family, implying the
result.

If $p\neq p^\Uparrow$, order elements of $\mathcal{HC}(p)$ by their
row number, and let $s$ be the greatest empty square of
$\mathcal{HC}_\mathfrak{D}$ preceding  the greatest filled square in
$\mathcal{HC}_\mathfrak{D}$. Let $t$ be the least filled square
following $s$ in $\mathcal{HC}_\mathfrak{D}$. The pair $\{s,t\}$ can
be extended to a cycle structure set for $p$ via Proposition
\ref{proposition:cs2}. Let $p^\uparrow$ be the partition obtained
from $p$ by filling $s$ and emptying $t$.  Then by Corollary \ref{corollary:module:tableau}, $p$ and $p^\uparrow$ lie in the same family. This process can be repeated successively
producing a sequence $p, p^\uparrow, (p^\uparrow)^\uparrow, \ldots$
of partitions in the same family which terminates in $p^\Uparrow.$
\end{proof}

\begin{example} Consider the symbol
$$
\Lambda = \left( {0\;1\; 3\; 4 \atop 2} \right)
$$
of defect $s=3$.  Its set of singles has four $3$-admissible
involutions $(0,1)$, $(1,2)$, $(2,3)$, and $(3,4)$ which, according to the
above proposition, produce the constructible representations
$$
\mathcal{S}_{(0,1)} : \left[\left( {1\;2\; 3\; 4 \atop 0}
\right)\right] \oplus \left[\left({0\;2\; 3\; 4 \atop
1}\right)\right] \hspace{.2in} \mathcal{S}_{(1,2)} : \left[\left(
{0\;2\; 3\; 4 \atop 1} \right)\right] \oplus \left[\left({0\;1\; 3\;
4 \atop 2}\right)\right]
$$
$$
\mathcal{S}_{(2,3)} : \left[\left( {0\;1\; 3\; 4 \atop 2}
\right)\right] \oplus \left[\left({0\;1\; 2\; 4 \atop
3}\right)\right] \hspace{.2in} \mathcal{S}_{(3,4)} : \left[\left(
{0\;1\; 2\; 4 \atop 3} \right)\right] \oplus \left[\left({0\;1\; 2\;
3 \atop 4}\right)\right]
$$
By using the identification from Section \ref{section:symbols}, we
can rephrase this list in terms of partitions of rank $r=s-1=2$. The
symbol $\Lambda$ corresponds to the partition $(4,3^2,1)$ and the
constructible representations can be rewritten in terms of
partitions as
$$
\mathcal{S}_{(0,1)}: [( 4,3,2,1^2 )] \oplus [(4,3,2^2
)]\hspace{.2in} \mathcal{S}_{(1,2)}: [( 4,3,2^2 )] \oplus [(4,3^2,1
)]
$$
$$
\mathcal{S}_{(2,3)}: \phantom{44} [( 4,3^2,1 )] \oplus [(4^2,2,1 )]
\hspace{.2in} \mathcal{S}_{(3,4)}: [( 4^2,2,1 )] \oplus [(5,3,2,1 )]
$$
\end{example}

\end{document}